\documentclass[a4paper,12pt]{article}
\usepackage{amsmath,amsfonts,amssymb, dsfont}
\usepackage{graphicx}

\newcommand {\fer}   {\eqref}

\newcommand {\da}   {\delta}

\newcommand {\f}   {\frac}
\newcommand {\p}   {\partial}

\newcommand{\beq}{\begin{equation}}
\newcommand{\eeq}{\end{equation}}

\newcommand{\qed}{{ \hfill
                     {\unskip\kern 6pt\penalty 500 \raise -2pt\hbox{\vrule\vbox to 6pt{\hrule width 6pt
                     \vfill\hrule}\vrule} \par}   }}

\newtheorem{lem}{Lemma}
\newtheorem{rem}{Remark}

\newtheorem{prop}{Proposition} 
\newtheorem{definition}{Definition} 
\newcommand{\R}{{\mathbb R}}

\begin{document} 
 
\title{Population structured by a space variable and a phenotypical trait}

\author{
Sepideh Mirrahimi\\
\footnotesize{UPMC, CNRS UMR 7598, Laboratoire Jacques-Louis Lions, F-75005, Paris.}\\
\footnotesize{ mirrahimi@ann.jussieu.fr}\\
\\
Ga\"el Raoul\\
\footnotesize{DAMTP,Centre for Mathematical Sciences, Cambridge University,}\\
\footnotesize{Wilberforce Road, Cambridge CB3 0WA, UK. }\\
\footnotesize{ g.raoul@damtp.cam.ac.uk}
}
\maketitle 
 
\begin{abstract}
We consider populations structured by a phenotypic trait and a space variable, in a non-homogeneous environment. In the case of sexual populations, we are able to derive models close to existing models in theoretical biology, from a structured population model. We then analyze the dynamics of the population using a simplified model, where the population either propagates through the whole space or it survives but remains confined in a limited range. For asexual populations, we show that the dynamics are simpler. In this case, the population cannot remain confined in a limited range, i. e. the population, if it does not get extinct, propagates through the whole space.
\end{abstract}

\section{Introduction}

In this paper, we are interested in populations that are structured by a continuous phenotypic trait $v\in \mathbb R$ and a continuous space variable $x\in\mathbb R$, living in a non-homogeneous environment: we will consider a phenotypic trait of best adaptation $\theta(x)$ that depends on the space variable. This type of population have been studied in \cite{PLB,KB,PYW}, using mostly numerical simulations.

\medskip

This type of models can for instance be used to study the distribution of a species along an environmental gradient (such as the north-south gradient of temperature in the northern hemisphere): To study the range of the species and its local adaptation, one should consider both evolution and spatial dynamics (see \cite{HGCBS,KB,PBM,BESP}). Those models are also useful to study the impact of an environmental change (e.g. global warming) on a population (see \cite{PLB,KB,PBM}).

\bigskip

Our work is largely based on \cite{KB}, and related articles \cite{PLB,KB,PYW,PBM}. In \cite{KB}, a partial differential equation model describing the spatial and evolutionary dynamics of a population has been introduced: 
\begin{equation}\label{eqbarton}
\left\{\begin{array}{l}
\partial_t N(t,x)-\Delta_xN(t,x)=\left(1-\frac 12(Z(t,x)-Bx)^2-N(t,x)\right)N(t,x),\\
\partial_tZ(t,x)-\Delta_x Z(t,x)=2\partial_x(\log N(t,x))\partial_xZ(t,x)+A(Bx-Z(t,x)).
\end{array}\right.
\end{equation}
where $N(t,x)$ is the population density at the location $x$, $Z(t,x)$ its mean phenotypic trait, and $A,\,B$ are two parameters. Numerical simulations where run for this model, and they showed that depending on $A$ and $B$, three biological scenarios were possible:
\begin{itemize}
\item if $B$ is large (the environment changes rapidly in space), then the population goes extinct,
\item for intermediate values of $B$, the population survives, but remains in a limited area,
\item if $B$ is small, the population invades the whole space.
\end{itemize}


\bigskip

In this paper, we show how \eqref{eqbarton} (indeed, the closely related model \eqref{eq2}) can be derived from a structured population model in the case of sexual populations, but isn't appropriate to model asexual populations. We also introduce a simplified model that allows us to investigate the dynamics of the population.

\bigskip

In Section 2 , we introduce a structured population model for the evolution of sexual populations structured by both a phenotypic trait and a space variable. To construct this model, we add a spatial variable to a well estabilshed local selection-mutation model (similarly to many kinetic models in physics or chemistry, see \cite{Villani}). This structured population model can also be seen as a continuous version of the model introduced in \cite{PYW}. In this structured population model, three parameters appear: $A,\,B,\,C$. We show that a model very close to \eqref{eqbarton} can be obtained as a formal limit  of our structured population model when $C$ is large (which means that many generations are necessary to obtain a significant growth of the population), see \eqref{eq2}. We explain in Rem. \ref{factorV} why we couldn't obtain exactly the same model.

Moreover, we consider a formal limit of \eqref{eq2} when $A$ is small (which corresponds to the values of $A$ considered in \cite{KB}). This way, we obtain a simpler model on $Z(t,x)$ only:
\begin{eqnarray}
\partial_tZ(t,x)-\Delta_x Z(t,x)&=&-4\frac{(\partial_x Z(t,x)-B/\sqrt{2})(Z(t,x)-(B/\sqrt{2})x)}{1-(Z(t,x)-(B/\sqrt{2})x)^2}\partial_xZ(t,x)\nonumber\\
&&+((B/\sqrt{2})x-Z(t,x)).\label{eq:modelZ}
\end{eqnarray}

\medskip

In Section 3, we analyse the model \eqref{eq:modelZ} derived in Section 2. Unfortunately, this equation has singularities that are obstacles to have a well-defined problem: we show that viscosity solutions exist, but are not unique. Nevertheless, this simple model allows us to describe two of the three possible scenarios from \cite{KB}: invasion fronts, and populations remaining in a limited area. The extinction phenomena cannot be observed here because of our assumption that $A$ is small.

\medskip

In Section 4, we investigate the case of asexual populations. Similarly as in Section 2, we introduce a structured population model for asexual populations structured by both a phenotypic trait and a space variable. After a rescaling, we obtain a structured population model that depends on two parameters, $A,\, B$. We show that in this case, only two biological scenarios are possible: either the population goes extinct, or it spreads to the entire space (we were however only able to show this last result in a weak sense, see Thm. \ref{extb}). This result shows in particular that the model \eqref{eqbarton} from \cite{KB} does not apply to asexual populations.

\section{Sexual populations}

\subsection{The model}

We start from a classical model describing the evolution of a population structured by a phenotypic trait only (see e.g. \cite{Burger,DBLD,MNG}, and \cite{DJMP,DJMR,LMP} for mathematics properties of this kind of models). Let $n(t,v)$ be the density of the population at time $t\geq 0$ and phenotypic trait $v\in \mathbb R$. We assume that the fitness depends on the square of distance between $v$ and an optimal adaptation trait $\theta$, and is altered by the population size. Then, the fitness $s[n(t,\cdot)](v)$ of an individual of phenotypic trait $v$ living among a resident population $n(t,\cdot)$ is then given by:
\begin{equation}\label{fitness}
s[n(t,\cdot)](v)=r_{max}-\frac 1{2V_s}(v-\theta)^2-\frac {r_{max}}K\int n(t,w)\,dw,
\end{equation}
where $r_{max}>0$ denotes the maximal growth rate of the population, $V_s$ measures the strength of the selection toward the optimal trait $\theta$, and $K$ is the carrying capacity of the environment.

\medskip

The fitness is well-defined in the case of asexual populations: it is the rate of births of offspring of trait $x$ minus the rate of death. In the case of sexual populations, however, the situation is more complicated, since reproduction requires two parents, of traits $v_\ast$ and $v_\ast'$, that give birth to an offspring of trait $v$, usually different from $v_\ast$ or $v_\ast'$. This has two consequences (see \cite{Bulmer}):
\begin{itemize}
\item We then need to define $Q(\cdot,v_\ast,v_\ast')$ the distribution function of the trait of the offspring. We will analyze the properties of $Q$ in Subsection \ref{subQ}.
\item We have to define the birth rate and the death rate separately. We will assume here that the birth rate is a constant, equal to $\gamma>r_{max}$.
\end{itemize}
If moreover we do not differentiate male and female (for instance because they have the same distribution), assume that mating is random and uniform among the population, and that the number of offspring is proportional to the population density (the idea being that the number of births is proportional to the number of females), then the evolution of the population structured by a phenotypic trait only is described by (see \cite{DBLD}):
\begin{eqnarray}
\partial_t n(t,v)&=&\left[-(\gamma-r_{max})-\frac {I_s}2-\frac 1{2V_s}(v-\theta)^2-\frac {r_{max}}K\int n(t,w)\,dw\right]n(t,v)\nonumber\\
&&+\gamma\int\int\frac {n(t,v_\ast)n(t,v_\ast')}{\int n(t,w)\,dw}Q(v,v_\ast,v_\ast')\,dv_\ast\,dv_\ast',\label{sexhomogene}
\end{eqnarray}
where $\gamma$ is the birth rate ($\gamma\geq r_{max}$), that we assume constant within the population (the selection occurs in the death term), and $\frac {I_s}2$ the additional death rate due to lethal mutations (see \cite{KB}).

\medskip

We consider next populations that are structured by a phenotypic trait $v$ as above, but also by a space variable $x\in \mathbb R$: $n(t,x,v)$. We assume that the selection-mutation process described above occurs locally in space, but that individuals move randomly in space (which we model by a diffusion of rate $\sigma_x$), and that the trait of optimal adaptation changes linearly in space:
\begin{equation}\label{b}
\theta(x)=bx.
\end{equation}
We then get the following model for sexual populations:
\begin{align}
&\partial_t n(t,x,v)-\sigma_x^2\Delta_xn(t,x,v)\nonumber\\
&\quad=\left[-(\gamma-r_{max})-\frac {I_s}2-\frac 1{2V_s}(v-bx)^2-\frac {r_{max}}K\int n(t,x,w)\,dw\right]n(t,x,v)\nonumber\\
&\qquad+\gamma\int\int\frac {n(t,x,v_\ast)n(t,x,v_\ast')}{\int n(t,x,w)\,dw}Q(v,v_\ast,v_\ast')\,dv_\ast\,dv_\ast'.\label{equnscaled}
\end{align}
For an existence theory for \eqref{equnscaled}, we refer to \cite{Prevost}.

\subsection{Properties of the sexual reproduction kernel $Q$}\label{subQ}

In Subsection \ref{resc}, we will rescale \eqref{equnscaled} to show that it indeed only depends on three parameters. But to be able to do so, we need first to define more precisely the reproduction kernel $Q$, and analyze its properties.


The sexual reproduction has two opposite effects on the repartition of the population over the phenotypic trait:
\begin{itemize}
\item For each set of two chromosomes, one come from the parent of phenotypic trait $v_\ast$, and one from the parent of trait $v_\ast'$, this process tends to give the offspring a trait between $v_\ast$ and $v_\ast'$. The effect of this phenomenon is to concentrate the population traits.
\item On the contrary, a variability is maintained in the population by mutations and recombinations (notice that the effect of recombinations is typically much larger than the effect of mutations, see \cite{Burger}).
\end{itemize}


A reasonable assumption is then to assume that in absence of selection, a sexual population phenotypic distribution will converge to a given profile that only depends on $Q$. To make this assumption precise, we consider the homogeneous sexual reproduction model \eqref{sexhomogene} without selection, and with a constant population size, that is
\begin{equation}\label{sexual}
\left\{\begin{array}{l}
\partial_t n(t,v)=\int\int Q(v,v_\ast,v_\ast')n(t,v_\ast)n(t,v_\ast')\,dv_\ast\,dv_\ast'-n(t,v)\\
n(0,v)=n^0(v)\in L^1(\mathbb R), \textrm{ with }\int n^0(v)\,dv=1,
\end{array}\right.
\end{equation}
then we assume that the long-time dynamics of this model is simple in the sense that:
\bigskip

\noindent {\bf Assumption 1:} \begin{itemize}
\item for any $v_\ast,\,v_\ast'\in \mathbb R$, $\int Q(v,v_\ast,v_\ast')\,dv=1$,
\item for any $v_\ast,\,v_\ast'\in \mathbb R$, $\int v\,Q(v,v_\ast,v_\ast')\,dv=\frac {v_\ast+v_\ast'}2$,
\item There exists $G$ (from the \emph{genetic variance}, see \cite{KB}) such that for any initial population $n^0$, the phenotypic variance of the population converges to $G$:
\begin{equation*}
\int \left(v-\int w\,n^0(w)\,dw\right)^2n(t,v)\,dv\to G^2,\textrm{ as }t\to\infty,
\end{equation*}
\item For any initial population $n^0$, the third moment of the phenotypic distribution of the population converges to 0:
\begin{equation*}
\int \left(v-\int w\,n^0(w)\,dw\right)^3n(t,v)\,dv\to 0,\textrm{ as }t\to\infty.
\end{equation*}
\end{itemize}

The reproduction kernel that is typically used  in structured population models for sexual populations (see \cite{Bulmer,DBLD}) is:
\begin{equation}\label{eqexQ}
Q(v,v_\ast,v_\ast')
:=\frac 1 {\gamma \sqrt{2\pi}}e^{-\frac {\left(v-\frac{v_\ast+v_\ast'}2\right)^2}{2\gamma^2}}.
\end{equation}
We show that Assumption 1 is satisfied for this particular reproduction kernel:
\begin{prop}\label{exQ}
Assumption 1 is satisfied by the reproduction kernel defined in \eqref{eqexQ}, with $G=\sqrt 2 \gamma$
\end{prop}

\begin{rem}
This assumption is indeed true under the more general assumption that $Q(v,v_\ast,v_\ast')=\left(\tilde \Gamma \ast \tilde Q(\cdot,v_\ast,v_\ast')\right)(v)$, where $\tilde \Gamma$ is symetrical with a positive variance, $\tilde Q(\frac{v_\ast+v_\ast'}2+v,v_\ast,v_\ast')=\tilde Q(\frac{v_\ast+v_\ast'}2-v,v_\ast,v_\ast')$ for any $v,\,v_\ast,\,v_\ast'\in \mathbb R$, 
and 
\begin{equation*}
\int \left(v-\frac {v_\ast+v_\ast'}2\right)^2\tilde Q(v,v_\ast,v_\ast')\,dv= C\frac {(v_\ast-v_\ast')^2}4,
\end{equation*}
with $C<1$. For this and more on this type of problem, we refer to \cite{HFMW,MRZ}.
\end{rem}

\bigskip

\noindent\textbf{Proof of Prop. \ref{exQ}:}

The two first conditions of Assumption 1 can be easilly checked. We then only prove the two last ones.

If we assume w.l.o.g. that $\int w\,n^0(w)=0$, then
\begin{eqnarray*}
\frac d{dt}\int v^2n(t,v)\,dv&=&\int\int \left(\int v^2Q(v,v_\ast,v_\ast')\,dv\right)n(t,v_\ast')n(t,v_\ast)\,dv_\ast\,dv_\ast'-\int v^2n(t,v)dv\\
&=&\int\int \left(\gamma^2+\left(\frac{v_\ast+v_\ast'}2\right)^2\right)n(t,v_\ast')n(t,v_\ast)\,dv_\ast\,dv_\ast'-\int v^2n(t,v)dv\\
&=&\gamma^2-\frac 12\int v^2n(t,v)dv.
\end{eqnarray*}
We deduce that the variance of $n$ converges to $G=\sqrt 2\gamma$ for any initial condition $n^0$:
\begin{equation*}
\int \left(v-\int w\,n^0(w)\,dw\right)n(t,v)\,dv\to G^2=2\gamma^2,\textrm{ as }t\to\infty.
\end{equation*}
Similarly we have
\begin{eqnarray*}
\frac d{dt}\int v^3n(t,v)\,dv&=&\int\int \left(\int v^3Q(v,v_\ast,v_\ast')\,dv\right)n(t,v_\ast')\,dv_\ast\,dv_\ast'-\int v^3n(t,v)dv\\
&=&\int\int \left(\frac{v_\ast+v_\ast'}2\right)^3n(t,v_\ast')\,dv_\ast\,dv_\ast'-\int v^3n(t,v)dv\\
&=&-\frac 34\int v^3n(t,v)dv,
\end{eqnarray*}
so that $\int v^3n(t,v)\to 0$ as $t\to\infty$. 
\begin{flushright}
$\square$
\end{flushright}

\subsection{Rescaling of the structured population model}\label{resc}
%
%

To simplify \eqref{equnscaled}, we perform the following rescaling:
\begin{equation*}
\tilde t:=\left(r_{max}-\frac {G^2}{2V_s}-\frac {I_s}2\right)^{-1}t,
\end{equation*}
\begin{equation*}
\tilde x=\sigma_x\left(r_{max}-\frac {G^2}{2V_s}-\frac {I_s}2\right)^{-1/2}x,
\end{equation*}
\begin{equation*}
\tilde n=\frac{r_{max}}{KG}\left(r_{max}-\frac {G^2}{2V_s}-\frac {I_s}2\right)^{-1}n, 
\end{equation*}
\begin{equation*}
 \tilde v=\frac v{G},
\end{equation*}
\begin{equation*}
\tilde Q(v+\frac{v_\ast+v'_\ast}2,v_\ast,v'_\ast):= G \,Q\left(Gv,Gv_\ast,Gv'_\ast\right),
\end{equation*}
where $G$ is given by Assumption 1. Then, \eqref{equnscaled} then becomes:
\begin{eqnarray}
&&\partial_t n(t,x,v)-\Delta_xn(t,x,v)\nonumber\\
&&\quad =-\left[(C-\frac A 2)+\frac A2(v-Bx)^2+\int n(t,x,w)\,dw\right]n(t,x,v)\nonumber\\
&&\qquad +(C+1)\int\int \frac{n(t,x,v_\ast)n(t,x,v'_\ast)}{\int n(t,x,w)\,dw}Q(v,v_\ast,v'_\ast)\,dv_\ast\,dv'_\ast,\label{eqscaled}
\end{eqnarray}
where
\begin{equation*}
A:=\frac {G^2}{V_s}\left(r_{max}-\frac {G^2}{2V_s}-\frac {I_s}2\right)^{-1},
\end{equation*}
\begin{equation*}
B:=\frac {b}{G\sigma_x}\left(r_{max}-\frac {G^2}{2V_s}-\frac {I_s}2\right)^{1/2},
\end{equation*}
\begin{equation*}
C:=\gamma\left(r_{max}-\frac {G^2}{2V_s}-\frac {I_s}2\right)^{-1}-1.
\end{equation*}
Moreover, the rescaled reproduction kernel satisfies the Assumption 1 with $G=1$.
\begin{rem}\label{factorV}
This is not the only possible rescaling for this equation. As we will see in Subsection  \ref{moments}, this particular scaling will allow us to obtain \eqref{eq2} as a formal limit of \eqref{eqscaled} when $C$ is large.
%
In \cite{KB}, the authors perform a renormalisation on \eqref{eqbarton} to get only two parameters $A$ and $B$. Those parameters are then defined as:
\begin{equation*}
A=\frac G{2V_s(r_{max}-I_s/2)},\quad B=\frac {\sigma_xb}{\sqrt{2V_s}(r_{max}-I_s/2)}.
\end{equation*}
The scaling they use wouldn't work here, because they assume that the scaling of the variable $v$ that they use does not modify $G$, the typical phenotypic variance of the phenotypic distribution of the population. With \eqref{equnscaled}, $G$ is necessarly affected by a scaling in the $v$ variable (see Assumption 1). This is also why we couldn't obtain exactly \eqref{eqbarton}, but the slightly different model \eqref{eq2}
\end{rem}

\subsection{Formal limit of the structured population model}\label{moments}

We denote by $N,\,Z,\,V$ the following moments of the distribution $n(t,x,\cdot)$:
\begin{equation*}
N(t,x):=\int n(t,x,v)\,dv,\quad Z(t,x):=\int v \frac{n(t,x,v)}{N(t,x)}\,dv,
\end{equation*}
\begin{equation*}
V(t,x):=\int (v-Z(t,x))^2\frac{n(t,x,v)}{N(t,x)}\,dv.
\end{equation*}

We show that $N$ and $Z$ satisfy the following unclosed equations:

\begin{prop} \label{propNZ}
If $n$ is a solution of \eqref{eqscaled}, then the moments of the phenotypic distribution of $n$ satisfy:
\begin{align}
&\partial_t N(t,x)-\Delta_xN(t,x)\nonumber\\
&\quad=\left[1+\frac A2(1-V(t,x))-\frac A2(Z(t,x)-Bx)^2-N(t,x)\right]N(t,x).\label{eqN}
\end{align}
\begin{align}
&\partial_tZ(t,x)-\Delta_x Z(t,x)\nonumber\\
&\quad=2\partial_x(\log N(t,x))\partial_xZ(t,x)+A(Bx-Z(t,x))V(t,x)\nonumber\\
&\qquad-\int(v-Z)^3\frac {n(t,x,v)} {N(t,x)} \,dv.\label{eqZ}
\end{align}
\end{prop}

%

\bigskip

\noindent\textbf{Proof of Prop. \ref{propNZ}:}

To get \eqref{eqN}, we integrate \eqref{eqscaled} along $v$:
\begin{align*}
&\partial_tN(t,x)-\Delta_x N(t,x)\\
&\quad =\int \partial_tn(t,x,v)-\Delta_x n(t,x,v)\,dv\\
&\quad = \left[1+\frac A2-\frac A2(Z(t,x)-Bx)^2-N(t,x)\right]N(t,x) -\frac A2\int (v-Bx)^2n(t,x,v)\, dv\\
&\quad =\left[1+\frac A2(1-V(t,x))-\frac A2(Z(t,x)-Bx)^2-N(t,x)\right]N(t,x).
\end{align*}

The second equation, \eqref{eqZ}, is obtained as follows:
\begin{align*}
&\partial_tZ(t,x)-\Delta_x Z(t,x)\\
&\quad=\partial_t\int v\frac {n(t,x,v)} {N(t,x)} \,dv-\Delta_x \int v\frac {n(t,x,v)} {N(t,x)}\,dv\\ 
&\quad=\int \frac v {N(t,x)} \left(\partial_tn(t,x,v)-\Delta_x n(t,x,v)\right)\,dv\\ 
&\qquad-\left(\partial_t N(t,x)-\Delta_x N(t,x)\right)\frac {Z(t,x)} {N(t,x)}\\ 
&\qquad+2\partial_x(\log N(t,x))\partial_xZ(t,x),
\end{align*} 
and then,
\begin{align*}
&\partial_tZ(t,x)-\Delta_x Z(t,x)\\
&\quad= -\int \frac v {N(t,x)} \left(\left(C-\frac A2\right)+\frac A2(v-Bx)^2+\int n(t,w) \,dw\right)n(t,x,v)\,dv\\
&\qquad+(C+1)\int \frac v{N(t,x)}\left(\int\int\frac{n(t,x,v_\ast)n(t,x,v_\ast')}{\int n(t,x,w)\,dw}\rho\ast_v Q(v,v_\ast,v_\ast')\,dv_\ast\, dv_\ast'\right)\,dv\\ 
&\qquad-\left(1-\frac A2(Z(t,x)-Bx)^2-N+\frac A2(1-V(t,x))\right)N(t,x)\frac {Z(t,x)} {N(t,x)} \\
&\qquad+2\partial_x(\log N(t,x))\partial_xZ(t,x)\\
&\quad= 2\partial_x(\log N(t,x))\partial_xZ(t,x)+A(Bx-Z(t,x))V(t,x)\\
&\qquad-\int(v-Z)^3\frac {n(t,x,v)} {N(t,x)} \,dv,
\end{align*} 
where we have used the fact that the reproduction kernel does not affect the mean phenotypic trait $\int v Q(v,v',v'_\ast)\,dv=\frac{v'+v'_\ast}2$.
\begin{flushright}
$\square$
\end{flushright}

\begin{rem}\label{geneflow}
The term $2\partial_x(\log N)\partial_xZ$ is referred to by biologists as the "gene flow" term (see \cite{Mayr} and \cite{PLB,KB,PBM}). This term models the fact that the mean phenotype of low density areas are greatly influenced by the phenotypes of neighboring high density areas. It is interesting to notice that this term does not come from the sexual reproduction term, but from the diffusion term:
\begin{equation*}
\Delta Z(t,x)=\int v\frac {\Delta_x n(t,x,v)}{N(t,x)}\,dv-\frac {Z(t,x)}{N(t,x)}\Delta N(t,x)+2\partial_x(\log N(t,x))\partial_xZ(t,x).
\end{equation*}
\end{rem}

To close the equations on $N$ and $Z$ obtained in Prop. \ref{propNZ}, notice that \eqref{eqscaled} can be written:
\begin{eqnarray*}
&&\partial_t n(t,x,v)-\Delta_xn(t,x,v)\\
&&\quad =C\left[\int\int \frac{n(t,x,v_\ast)n(t,x,v'_\ast)}{\int n(t,x,w)\,dw}Q(v,v_\ast,v'_\ast)\,dv_\ast\,dv'_\ast-n(t,x,v)\right]\\
&&\qquad +\left[\frac A2-\frac A2(v-Bx)^2-\int n(t,x,w)\,dw\right]n(t,x,v)\\
&&\qquad +\int\int \frac{n(t,x,v_\ast)n(t,x,v'_\ast)}{\int n(t,x,w)\,dw}Q(v,v_\ast,v'_\ast)\,dv_\ast\,dv'_\ast
\end{eqnarray*}
so that if $C$ is very large, the first term will dominate the dynamics of the population. Since this first term corresponds to the "pure" sexual reproduction equation \eqref{sexual}, and thanks to Assumption 1, it is natural to assume that at all time $t>0$ and all locations $x\in \mathbb R$,
\begin{align}
&V(t,x)=\int (v-Z(t,x))^2\frac{n(t,x,v)}{N(t,x)}\,dv\sim G=1,\nonumber\\
& \int (v-Z(t,x))^3\frac{n(t,x,v)}{N(t,x)}\,dv\sim 0. \label{closure}
\end{align}
Notice that here, thanks to the rescaling performed in the preceding subsection, $G=1$.


\medskip

 If we use those properties to close the system of equations on $N$ and $Z$, we get the model:
\begin{equation}\label{eq2}
\left\{\begin{array}{l}
\partial_t N(t,x)-\Delta_xN(t,x)=\left(1-\frac A2(Z(t,x)-Bx)^2-N(t,x)\right)N(t,x),\\
\partial_tZ(t,x)-\Delta_x Z(t,x)=2\partial_x(\log N(t,x))\partial_xZ(t,x)+A(Bx-Z(t,x)).
\end{array}\right.
\end{equation}
This model is very close to the model \eqref{eqbarton} from \cite{PLB,KB}. The model \eqref{eqbarton} was build directly, without the intermediate step of a structured population model, and the limits of its applications was unclear (see \cite{PYW}). Our derivation shows that the model \eqref{eq2} is valid (in the sense that it is the formal limit of \eqref{equnscaled}) if:
\begin{itemize}
\item The reproduction is sexual,
\item The reproduction kernel satisfies Assumption 1,
\item $C$ is large.
\end{itemize} 

\begin{rem}
$C=\gamma\left(r_{max}-\frac {G^2}{2V_s}-\frac {I_s}2\right)^{-1}-1$ is large if the birth rate is large compared to the maximal fitness of the population (or many generations are necessary to obtain a significant growth of the population, which seems reasonable in many biological situations). 

Notice also that for  the example of Prop. \ref{exQ}, the convergence of Assumption 1 are exponentially fast, so that the simplification \eqref{closure} may be accurate even if $C$ is not very large.
\end{rem} 

The model \eqref{eq2} may hold in other situations, justifying the closure assumption \eqref{closure} with other arguments. However, we show in Section \ref{section:as} that \eqref{eq2} cannot hold for asexual populations for the whole range of parameters $A$ and $B$ (see Remark \ref{eq2as}).

\subsection{Derivation of a simplified model}

To simplify \eqref{eq2}, we first apply the following change of variable: 
\begin{equation*}
\tilde t:=At,\quad \tilde x=\sqrt A x,
\end{equation*}
\begin{equation*}
\tilde Z=\sqrt{\f{A}{2}}Z.
\end{equation*}
The rescaled model becomes
\begin{equation*}
\left\{\begin{array}{l}
\partial_t N(t,x)-\Delta_xN(t,x)=\frac 1 { A}\left(1-(Z(t,x)-(B/\sqrt 2)x)^2-N(t,x)\right)N(t,x),\\
\partial_tZ(t,x)-\Delta_x Z(t,x)=2\partial_x(\log N(t,x))\partial_xZ(t,x)+((B/\sqrt{2})x-Z(t,x)).
\end{array}\right.
\end{equation*}
Now, if we assume that $A$ is very small, $N$ and $Z$ are related by the simple relation:
\begin{equation}\label{N}
N(t,x)\sim 1-(Z(t,x)-(B/\sqrt{2}) x)^2.
\end{equation}
Therefore we get the simpler model \eqref{eq:modelZ}, on $Z$ only:
\begin{eqnarray*}
\partial_tZ(t,x)-\Delta_x Z(t,x)&=&-4\frac{(\partial_x Z(t,x)-B/\sqrt{2})(Z(t,x)-(B/\sqrt{2})x)}{1-(Z(t,x)-(B/\sqrt{2})x)^2}\partial_xZ(t,x)\nonumber\\
&&+((B/\sqrt{2})x-Z(t,x)).
\end{eqnarray*}


\begin{rem}
In \cite{KB}, the range of $A$ that has been considered was $A\in [0.001,1]$. Our approximation assuming that $A$ is small thus seems reasonable.

\medskip

Another simplification  had been proposed in \cite{KB}, where the equation on $N$ was replaced by 
\begin{equation}\label{nexp}
N:=k\exp\left(\gamma \left(1-A(Z-Bx)^2\right)\right).
\end{equation}
\eqref{eqbarton} then simplifies considerably:
\begin{equation*}
\partial_t Z(t,x)-\Delta_x Z(t,x)=A(Bx-Z(t,x))\left[1-4\gamma\partial_xZ(t,x)(B-\partial_x Z(t,x))\right].
\end{equation*}
However, the simplification \eqref{nexp} seems independent of \eqref{eqbarton}. Our simplification has the advantage to rely on a clearer assumption: \eqref{eq:modelW} is the formal limit of \eqref{eq2} when $A$ is small. 
\end{rem}

\section{Dynamics of sexual populations}

\subsection{Well-posedness of the model}

By replacing $W=Z-(B/\sqrt{2})x$ in \fer{eq:modelZ} we obtain the following equation
\beq
\label{eq:modelW}
\p_t W-\triangle_x W=-4\f{\p_x W W }{1-W^2}(\p_xW+B/\sqrt{2})-W,
\eeq
with $-1\leq W\leq 1$. This equation has a singularity for $W=\pm 1$.
  The existence of singularities is an obstacle to have a well-defined problem. However, as we will see in section \ref{fronts}, the singularities lead to the existence of propagative fronts. In most of the cases in the classical study of propagative fronts, one proves the existence of propagative fronts that connect  an unstable steady state to a stable steady state. Here the situation is different. The propagative fronts connect the unstable steady state $W=0$ to the singular point $W=-1$. While the presence of singularities is crucial to observe propagative fronts, it is an obstacle to prove uniqueness or comparison results to compare the solution with the propagative fronts.  Nevertheless we are able to introduce an approached model where the uniqueness and comparison principles hold.

 Since \fer{eq:modelW} is singular, we approximate it by the following model
\beq
\label{approxmod}
\p_t W_\da-\triangle_x W_\da=-4\f{\p_x W_\da W_\da }{1-W_\da^2+\da}(\p_xW_\da+B/\sqrt 2)-\f{(1-W_\da^2)W_\da}{1-W_\da^2+\da},\\
\eeq
with
\beq \label{boundary}
W_\da(t=0,\cdot)=W^{0}_\da(\cdot).
\eeq
With this choice of approximation we avoid the singularities and transform the singularity in $-1$ into a stable steady state (the stability is for the ODE formulation presented in section \ref{fronts}).

Under the assumption
\beq
\label{w0bound}
- 1\leq W^{0}_\da\leq 1, 
\eeq
 equation \fer{approxmod} has a smooth solution that stays between $-1$ and $1$ by the maximum principle. 
Using the following assumption on the initial data
\beq\label{ass:lip}
|\p_x W_\da^0|\leq C_1,
\eeq
with $C_1$ a positive constant, we prove a uniform  Lipschitz bound for the $W_\da$'s and we deduce that the $W^\da$'s converge to a viscosity solution  of a variant of equation \fer{eq:modelW} (see \cite{C.I.L:92,GB:94} for general introduction to the theory of  viscosity solutions):
\begin{prop}\label{th:lip}
Under assumptions \fer{w0bound} and \fer{ass:lip}, we have that the solutions of \fer{approxmod} are uniformly bounded and Lipschitz:
$$-1\leq W_\da(t,x)\leq 1,\qquad |\p_x W_\da(t,x)|\leq C_2,\qquad\text{for all }(t,x)\in \R^+\times \R.$$
Consequently, after extraction of a subsequence, the $W_\da$'s converge to a continuous function $W$ that is a viscosity solution of 
\beq
\label{eq:modelW2}
(1-W^2)\p_t W-(1-W^2)\triangle_x W=-4{\p_x W W }(\p_xW+B/\sqrt 2)-(1-W^2)W.
\eeq
\end{prop}
We notice that equation \fer{eq:modelW2} is the original model \fer{eq:modelW} multiplied by $1-W^2$.\\

\bigskip

\noindent\textbf{Proof of Prop. \ref{th:lip}:}

 We differentiate equation \fer{eq:modelW} with respect to $x$ and obtain
\beq\label{pxW}\begin{array}{rl}
   \p_t \p_x W_\da-\triangle_x \p_x W_\da&=-4\f{ W_\da }{1-W_\da^2+\da}(2\p_xW_\da+B/\sqrt 2)\p_x(\p_x W_\da)\\
&-4\p_xW_\da^2(\p_x W_\da+B/\sqrt 2)\f{1+W_\da^2+\da}{(1-W_\da^2+\da)^2}\\
&-\p_xW_{\da}\left(1-\da\f{1+W_\da^2+\da}{(1-W_\da^2+\da)^2}\right),
  \end{array}
\eeq
where the last term comes from
$${\p_x} \left(\f{W_\da(1-W_\da^2)}{1-W_\da^2+\da}\right)={\p_x} \left(W_\da-\f{\da W_\da}{1-W_\da^2+\da}\right)=\p_xW_\da\left(1-\da\f{1-W_\da^2+\da+2W_\da^2}{(1-W_\da^2+\da)^2}\right).$$

We multiply \eqref{pxW} by $\p_x W_{\da}$ and devide by $|\p_xW_{\da}|$ and obtain
\beq\label{pxWda}\begin{array}{rl}
   \p_t |\p_xW_{\da}|-\triangle_x |\p_xW_{\da}|&\leq-4\f{ W_\da }{1-W_\da^2+\da}(2\p_xW_\da+B/\sqrt 2)\p_x(|\p_x W_\da|)\\
&-4\p_xW_\da^2(\p_x W_\da+B/\sqrt 2)\left(\f{1+W_\da^2+\da}{(1-W_\da^2+\da)^2}\right)\text{sgn}(\p_x W_\da)\\
&-\left(1-\da\f{1+W_\da^2+\da}{(1-W_\da^2+\da)^2}\right)
|\p_xW_{\da}|.
  \end{array}
\eeq
It follows that, for $\da<1$, $|\p_x W_\da|$ is a subsolution of the following equation
 \beq\label{eq:g}\begin{array}{rl}\p_t g-\triangle_x g&=\alpha(t,x)\p_x g+\left(\f{1+W_\da^2+\da}{(1-W_\da^2+\da)^2}\right)\left(-4 g^3+2\sqrt 2\,Bg^2+ g\right),\end{array}\eeq
 with $$\alpha(t,x)= -4\f{ W_\da }{1-W_\da^2+\da}\,(2\p_xW_\da+B/\sqrt 2). $$
 We choose a positive constant $C_2$ such that $-4 C_2^3+\sqrt 2\,BC_2^2+C_2<0$ and $C_1<C_2$. The constant $C_2$ is a supersolution of the equation above and thus
 $$|\p_xW_{\da}|\leq C_2.$$
We proved that the $W_\da$'s are uniformly Lipschitz continuous in space. Moreover we know that the $W_\da$'s are bounded. It follows that the $W_\da$'s are uniformly continuous in time (see \cite{GB.SB.OL:02}). Using Arzela Ascoli Theorem we conclude that, after extraction of a subsequence, the $W_\da$'s converge to a continuous function $W$. By the stability of viscosity solutions (see \cite{GB:94}), $W$ is a viscosity solution of \fer{eq:modelW2}.

\begin{flushright}
$\square$
\end{flushright}

\begin{rem} We can relax assumption \fer{ass:lip} in Proposition \ref{th:lip}. This is because 
$y=\f{1}{2\sqrt{t}}+L$ is also a supersolution of equation \fer{eq:g}, for $L=L(B)$ a large constant and $\da$ small. Therefore we have
 $$|\p_x W_\da|\leq \f{1}{2\sqrt{t}}+L.$$
It follows that there is a  regularizing effect and the $W_\da$'s become uniformly Lipschitz, for all $t>t_0>0$, even if they are not uniformly Lipschitz initially.
\end{rem}

We proved that the equation \fer{eq:modelW2} has a solution in the viscosity sense. Unfortunately the viscosity criteria is not enough to impose uniqueness. We give a counter-example below:\\

\noindent{\bf Example. Non-uniqueness for equation \fer{eq:modelW2}:}
We have the two following solutions to equation \fer{eq:modelW2}:
$$W_1(t,x)=-1,\qquad \text{for all }(t,x)\in \R^+\times \R,$$
$$W_2(t,x)=- e^{-t},\qquad \text{for all }(t,x)\in \R^+\times \R.$$
Here the biological solution is the first one. Because $W=-1 $ corresponds to $N=0$. Therefore, if initially $W(0,x)=-1 $ for all $x\in \R$, we expect that $W(t,\cdot)\equiv -1 $, for all $t\in R^+$. Otherwise some mass is created out of nowhere. 
We can easily verify that if $W_\da(t=0,\cdot)\equiv -1 $, we have $W_\da(t,\cdot)=-1 $ for all $t>0$. Therefore our approximation chooses the biological solution.\\

In section \ref{fronts} we study the propagative fronts for this model. To be able to compare the solutions with the propagative fronts and to show the propagation of the density in space, we need a comparison principle. Unfortunately as we saw above, the equation \fer{eq:modelW2} does not have a unique viscosity solution and therefore it does not admit a comparison principle. However we can prove a comparison principle for the approached model. We first recall it's definition:

\begin{definition}{ Comparison principle:} Equation $L(D^2 u, Du,u,x,t)=0$ verifies a comparison principle, if for any subsolution $w^1$ and supersolution $w^2$ of $L$ such that $w^1(0,x)\leq w^2(0,x)$, we have
$$w^1(t,x)\leq w^2(t,x),\qquad \text{for all }(t,x)\in\R^+\times \R.$$
\end{definition}

We prove that there is a comparison property for \fer{approxmod}. In particular \fer{approxmod} has a unique solution. 
\begin{prop}\label{th:comp}
 The problem \fer{approxmod} admits a comparison principle in the set of solutions $\{-1\leq W\leq 1\}$. 
\end{prop}

\bigskip

\noindent\textbf{Proof of Prop. \ref{th:comp}:}

We suppose that $W_1$ and $W_2$ are respectively subsolution and supersolution of \fer{approxmod} and  
$$W_1(t=0,\cdot)\leq W_2(t=0,\cdot).$$
We prove that $W_1\leq W_2$ for all $(t,x)\in \R^+\times \R$. Let $(\bar t,\bar x)$ a maximum point of $W_1-W_2$. Since it is a maximum point we
have $\p_xW_1(\bar t,\bar x)=\p_xW_2(\bar t,\bar x)=p$. Therefore we have
\begin{align*}
& \p_t (W_1-W_2)(\bar t,\bar x)-\triangle(W_1-W_2)(\bar t,\bar x)\\
&\quad \leq -4p\,(p+B/\sqrt 2)\f{(1+W_1W_2+\da)}{(1-W_1^2+\da)(1-W_2^2+\da)}\,(W_1-W_2)(\bar t,\bar x)\\
&\qquad-\f{(1+\da)(1-W_1^2-W_1W_2-W_2^2)+W_1W_2+W_1^2W_2^2}{(1-W_1^2+\da)(1-W_2^2+\da)}\,(W_1-W_2)(\bar t,\bar x).
  \end{align*}
In the previous section we proved that $|\p_x W|$ is bounded. Thus $p(p+B/\sqrt{2})$ is bounded. Moreover $W_1$ and $W_2$ are bounded and 
$$1-W_i^2+\da\geq \da,\qquad \text{for }i=1,\, 2.$$ Therefore the coefficient of $W_1-W_2$ is bounded. Following the classical maximum principle we deduce  that equation \fer{approxmod} admits a comparison principle.
\begin{flushright}
$\square$
\end{flushright}

\subsection{Existence of propagative fronts and steady populations}\label{fronts}

We are interested in propagation fronts, that is solutions of \eqref{eq:modelZ} of the form $Z(t,x)=(B/\sqrt 2)x+U(x-\nu t)$. The equation \eqref{eq:modelZ} becomes:
\begin{equation*}
-\nu U'-U''=-4\frac {U'U}{1-U^2}(U'+B/\sqrt{2})-U.
\end{equation*}
If we denote by $V:=U'$, finding a propagative front is then equivalent to find a solution defined on $\mathbb R$ to the ODE given by the vector field
\begin{equation}\label{FUV}
\left\{\begin{array}{l}
F_U(U,V)=V,\\
F_V(U,V)=-\nu V+4\frac {UV}{1-U^2}(V+B/ \sqrt{2})+U.
\end{array}\right.
\end{equation}

To have a meaning with respect to \eqref{eq:modelW}, those solutions must satisfy $u(t)\in [-1,1]$. 

\begin{prop}\label{EDOfronts}
For any $B>0$, there exists $\nu_B\in\mathbb R$ such that \eqref{eq:modelW} has a propagative front of speed $\nu$, $Z(t,x)=(B/\sqrt 2)x+U(x-\nu t)$, satisfying  
\begin{equation*}
U(x)\to 0\textrm{ as }x\to -\infty,\quad U(x)\to -1\textrm{ as }x\to +\infty,
\end{equation*}
if and only if $\nu>\nu_B$.

The propagative front with speed $\nu$ is unique (up to a translation), and $\nu_B$ is a decreasing function of $B$. 
\end{prop}

\begin{figure}
\label{fig:energy2decay}
\includegraphics[angle=0,scale=0.45]{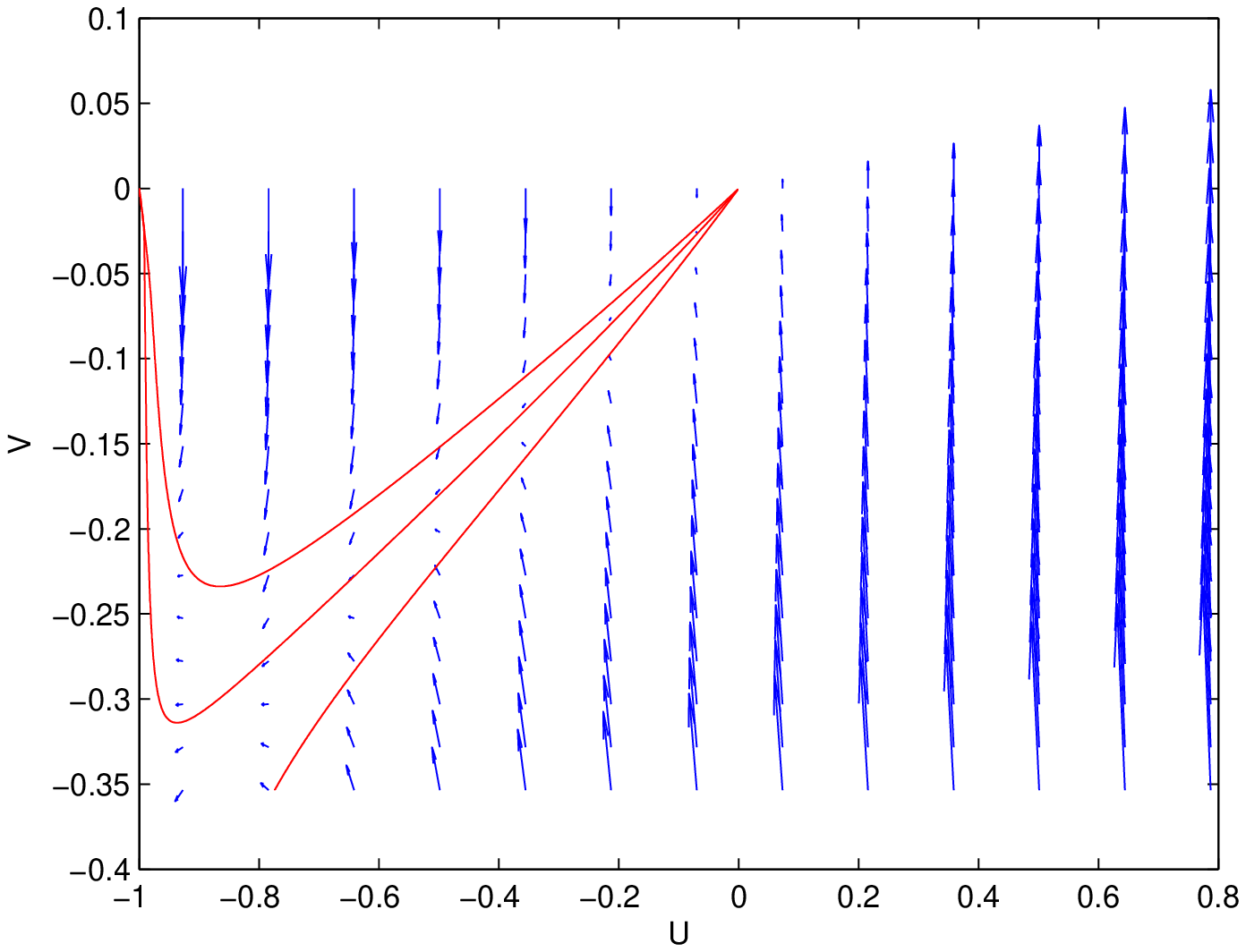}\includegraphics[angle=0,scale=0.45]{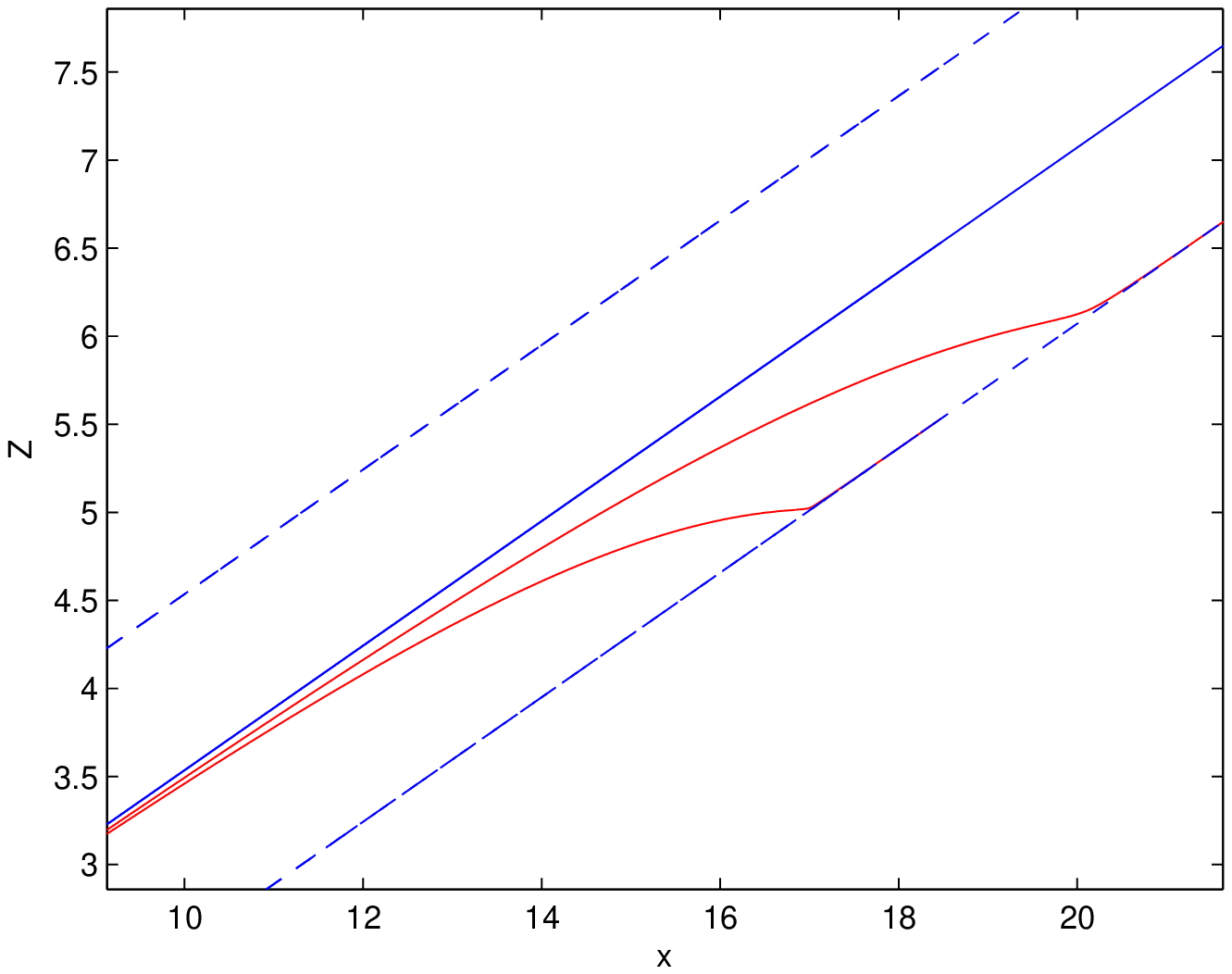}
\caption{On the left, we represent solutions of the ODE defined by the vector field \eqref{FUV} for $B:=0.5$ and $\nu=3.75,\,4.75,\,5.75$ (the vector field represented corresponds to $\nu=4.75$). A propagative front exists for those two last values of $\nu$ only, which we represent on the right.}
\end{figure}

For each $B>0$, there exists a one-parameter family of propagative fronts.
By analogy to the KPP-Fisher equation, one can guess that there is only one stable propagative front, the one with the least speed. Those fronts would be invasive fronts if $\nu_B>0$, and extinction fronts if $\nu_B<0$. The proposition \ref{EDOstat} shows that steady populations exist in this second case only.

We notice that, since the model does not admit a comparison principle, we cannot use the usual methods used in the study of the KPP-Fisher equation, to study the stability of propagative fronts rigorously. Nevertheless, if we choose those solutions of \fer{eq:modelW} that are obtained as the limits of the approached solutions $W_\da$'s,  and since the appraoched model admits a comparison principle, one can expect that the comparison principle be true for these limit solutions. The comparison principle would in particular allow us  to compare the solutions with the propagative fronts and prove the propagation of the population by the minimal spead of propagative fronts.

\begin{prop}\label{EDOstat}
The equation \eqref{eq:modelZ} has a non-trivial steady-state if and only if $\nu_B<0$. Steady-states $Z(t,x)=(B/\sqrt 2)x+U(x)$ satisfy 
\begin{equation*}
U(x)\to 1\textrm{ as }x\to -\infty,\quad U(x)\to -1\textrm{ as }x\to +\infty.
\end{equation*}
\end{prop}

\begin{figure}
\label{fig:energy2decay}
\includegraphics[angle=0,scale=0.45]{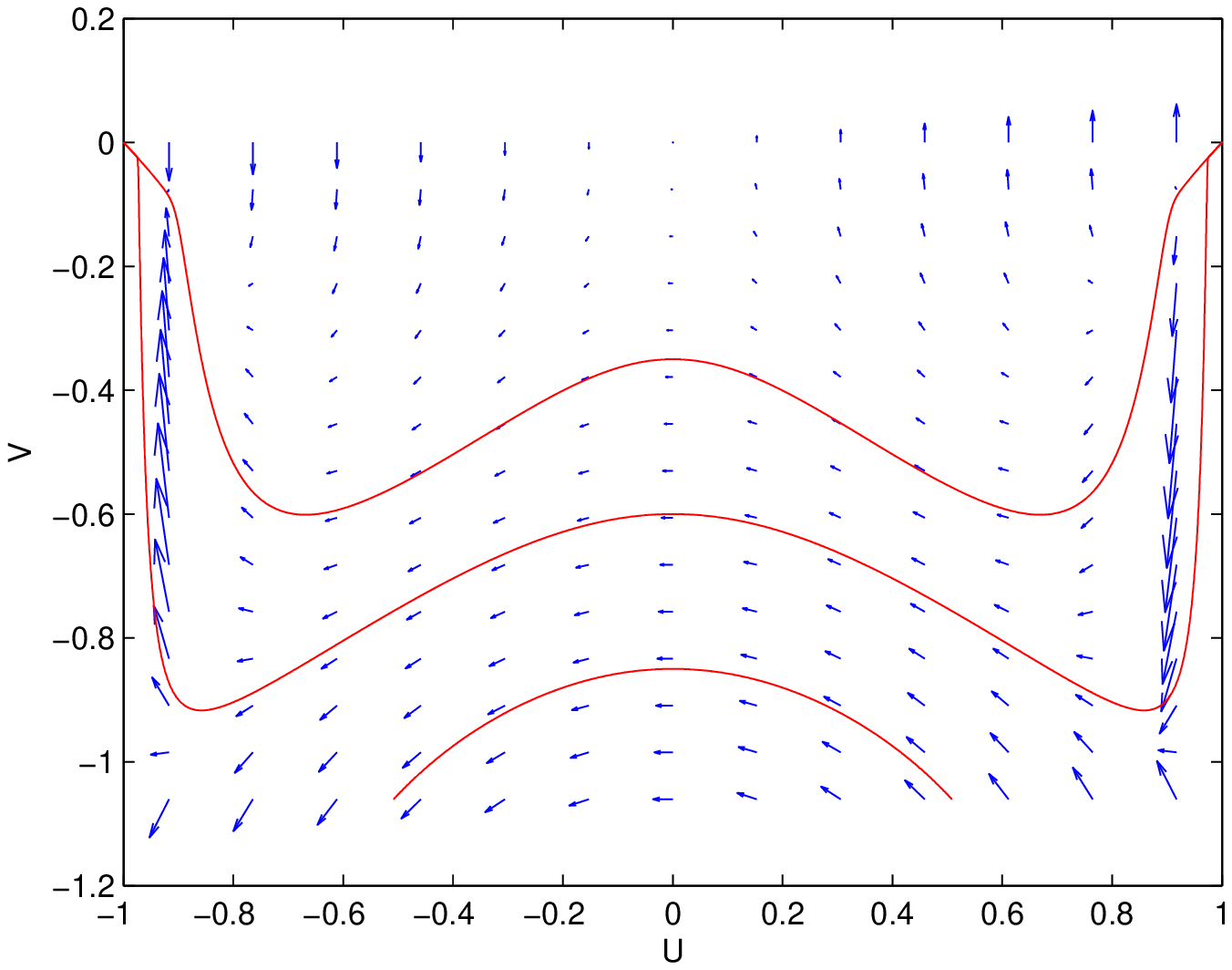}\includegraphics[angle=0,scale=0.45]{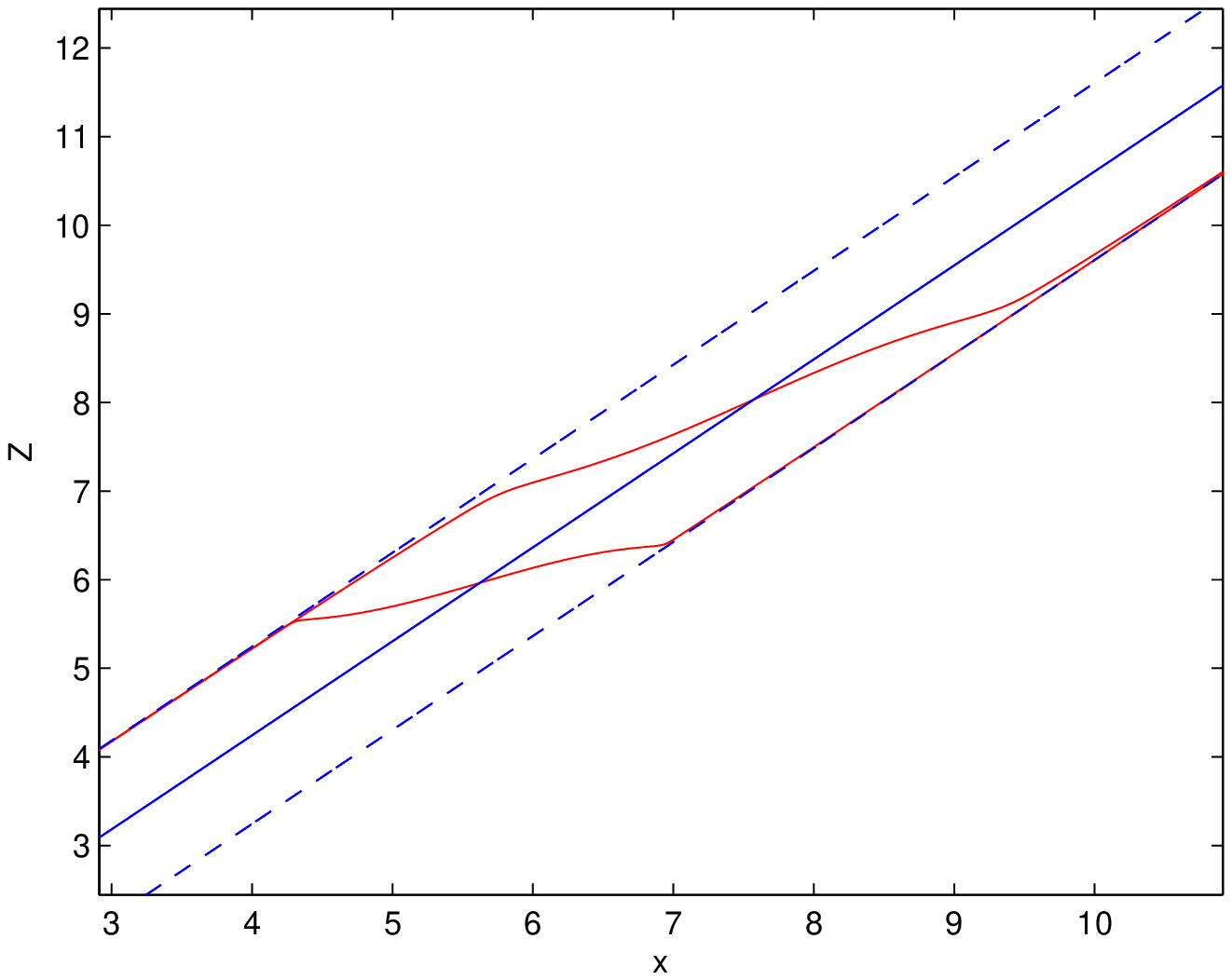}
\caption{On the left, we represent solutions of the ODE defined by the vector field \eqref{FUV} for $B:=1.5$ and $\nu=0$. Two of the three solutions represented define steady-states of \eqref{eq:modelW},  which we represent on the right. }
\end{figure}

As one can see in the proof, if $\nu_B<0$, there exists indeed a whole family of steady-states. If we assume that $U(0)=0$ (to avoid the translation invariance of the problem), then the family of steady-states can be parametrized by $U'(0)\in (-K_B,0)$, for some $K_B>0$.
%

To show those two propositions, we will use the two following lemma:

\begin{lem}\label{l1}
Let $\nu\in \mathbb R$. There exist two (up to a shift in the time variable) solutions $(u,v)$ to the ODE defined by the vector field \eqref{FUV} such that $(u,v)(t)\to (0,0)$ as $t\to -\infty$. At most one of them is globally defined, which satisfies:
\begin{equation*}
(u,v)(t)\sim - C_- e^{\frac {\sqrt{\nu^2+4}-\nu}2t}(1,\frac {\sqrt{\nu^2+4}-\nu}2).
\end{equation*}
For this solution, $u$ is strictly decreasing.
\end{lem}

\begin{lem}\label{l2}
Let $\nu_1\leq \nu_2$, and $(u_{\nu_1},v_{\nu_1}),\, (u_{\nu_1},v_{\nu_1})$ be the corresponding solutions given by Lemma \ref{l1}. If for $t_1,\,t_2\in \mathbb R$,
\begin{equation}\label{as:t1t2}\left\{
\begin{array}{l}
u_{\nu_1}(t_1)=u_{\nu_2}(t_2)\\
v_{\nu_1}(t_1)\leq v_{\nu_2}(t_2),
\end{array}\right.
\end{equation}
then, for any $t_1'> t_1$, $t_2'> t_2$ such that $u_{\nu_1}(t_1')=u_{\nu_2}(t_2')$, we have $v_{\nu_1}(t_1')\leq v_{\nu_2}(t_2')$, and this inequality is strict if $\nu_1< \nu_2$.
\end{lem}

\bigskip

\noindent\textbf{Proof of Lem. \ref{l1}:}

The Differential of the vector field $F$ in $(0,0)$ is
\begin{equation*}
DF_{(0,0)}=\left(\begin{array}{cc}
0&1\\
1&-\nu
\end{array}\right),
\end{equation*}
and $(0,0)$ is hyperbolic ($\det DF_{(0,0)}<0$). The Hartman-Grobman Theorem then applies, and there exists only two (non-trivial) solutions $(u,v)$ satisfying $(u,v)(t)\to_{t\to -\infty} (0,0)$. Since the eigenvector associated to the positive eigenvalue of $DF_{(0,0)}$ is $(1,\frac {\sqrt{\nu^2+4}-\nu}2)$, those two solutions are equivalent to
\begin{equation*}
(u,v)(t)\sim_{t\to -\infty} \pm C_\pm e^{\frac {\sqrt{\nu^2+4}-\nu}2t}\left(1,\frac {\sqrt{\nu^2+4}-\nu}2\right),
\end{equation*}
for some $C_-,\,C_+>0$.

\medskip

The solution such that $(u,v)(t)\sim_{t\to -\infty}  C_+e^{\frac {\sqrt{\nu^2+4}-\nu}2t}\left(1,\frac {\sqrt{\nu^2+4}-\nu}2\right)$ satisfies $u(\bar t)>0$, $v(\bar t)>0$ for some $\bar t$. Since $F_U(u(\bar t),V)>0$ for $V\geq v(\bar t)$ and $F_V(U,v(\bar t))\geq F_V(u(\bar t),v(\bar t))>0$ for $U\in [u(\bar t),1)$, the solution cannot escape $[u(\bar t),1)\times [v(\bar t),\infty)$. In particular, for $t\geq \bar t$, $u'(t)=v(t)\geq v(\bar t)>0$ and since the vector field is not defined for $U=1$, the solution cannot be global.

\medskip

The other solution satisfies $u(t)<0$, $v( t)<0$ for $t\geq \bar t$. Since $F_U(0,V)\leq0$ for $V\leq 0$ and $F_V(U,0)\leq 0$ for $U\leq 0$, the solution cannot escape $\mathbb R_-^2$, and in particular, $u'(t)=v(t)\leq 0$, which shows that $u$ is strictly decreasing at all times.

\begin{flushright}
$\square$
\end{flushright}

\bigskip

\noindent\textbf{Proof of Lem. \ref{l2}:}

We know that $u_\nu$ is strictly decreasing, we can thus define the graph of $(u_{\nu_1},v_{\nu_1})$.

We assume that $\bar t_1,\,\bar t_2$ are the smallest points respectively in $(t_1,\infty)$ and in $(t_2,\infty)$ such that $(u_{\nu_2},v_{\nu_2})(\bar t_2)=(u_{\nu_1},v_{\nu_1})(\bar t_1)$. We have,
\begin{eqnarray*}
v_{\nu_1}'(\bar t_1)&=&-\nu_1 v_{\nu_1}(\bar t_1)+4\frac {u_{\nu_1}(\bar t_1)v_{\nu_1}(\bar t_1)}{1-u_{\nu_1}(\bar t_1)^2}\left(v_{\nu_1}(\bar t_1)+B/\sqrt 2\right)+u_{\nu_1}(\bar t_1)\\
&=&v_{\nu_2}'(\bar t_2)+(\nu_2-\nu_1) v_{\nu_1}(\bar t_1)\\
&\leq&v_{\nu_2}'(\bar t_2),
\end{eqnarray*}
this inequality being strict if $\nu_1< \nu_2$. It follows that $\frac {v_{\nu_1}'(\bar t_1)}{u_{\nu_1}'(\bar t_1)}>\frac{v_{\nu_2}'(\bar t_2)}{u_{\nu_2}'(\bar t_2)}$. The graph of $(u_{\nu_2},v_{\nu_2})$ can thus only cross the graph of $(u_{\nu_1},v_{\nu_1})$ from below to above, when $t$ increases. This is enough to conclude that $v_{\nu_1}(t_1')\leq v_{\nu_2}(t_2')$, for all $t'_1>t_1$ and $t'_2>t_2$. Moreover, the latter inequality is strict if $\nu_1< \nu_2$. This completes the proof of Lem. \ref{l2}.
%
%
%
%
%
%
%
\begin{flushright}
$\square$
\end{flushright}

\bigskip

\noindent\textbf{Proof of Prop. \ref{EDOfronts}:}

\emph{Step 1:} We show that the solution $(u,v)$ given by Lem. \ref{l1} satisfies either $v(\bar t)=-B/\sqrt{2}$ for some $\bar t$, or $(u,v)(t)\to (-1,0)$ as $t\to +\infty$. Moreover, the solution is global in this last case only.

\bigskip

Since $u$ is strictly decreasing and $F_V(-1,V)=\infty$ for $V\in (-B/\sqrt 2,0)$, only two situations are possible: either $v(\bar t)=-B/\sqrt 2$ for some $\bar t<\infty$, or $(u,v)(t)\to (-1,0)$ as $t\to \bar t\in \mathbb R\cup\{\infty\}$.

\medskip

If $v(\bar t)=-B/\sqrt 2$, then $v(t)<-B/\sqrt 2$ for all $t\geq \bar t$. This is because $F_V(U,-B/\sqrt 2)=\nu B/\sqrt 2+U\leq F_V(u(\bar t),-B/\sqrt 2)<0$ for $U\in (-1,u(\bar t))$ and since $u$ is decreasing. Therefore, $u'(t)=v(t)<-B/\sqrt 2$ for all $t\geq \bar t$. From the latter together with $F_V(-1,V)=-\infty$ for $V<-B/\sqrt 2$ we obtain that the solution cannot be globally defined.

\medskip

Let $(u,v)(t)\to (-1,0)$ as $t\to \bar t\in \mathbb R\cup\{\infty\}$. We show that $\bar t=\infty$. For $(U,V)$ close to $(-1,0)$,  we have $F_V(U,V)\sim \frac {-2B}{\sqrt 2}\frac V{1+U}-1$. Then, 
\begin{eqnarray*}
\frac d {dt} \left(\frac v{1+u}\right)(t)&=&\frac {F_V(u(t),v(t))(1+u(t))-v(t)F_U(u(t),v(t))}{(1+u(t))^2}\\
&\sim&\frac 1{(1+u(t))^2}\left[\left(\frac {-2B}{\sqrt 2}\frac {v(t)}{1+u(t)}-1\right)(1+u(t))-v(t)^2\right]\\
&\sim& \frac 1{(1+u(t))^2}\left[\frac {-2B}{\sqrt 2}v(t)-v(t)^2-(1+u(t))\right]\\
&\geq& 0,
\end{eqnarray*}
if $v(t)\leq \frac {-\sqrt 2}{3B}(1+u(t))$ and $(u(t),v(t))$ is close to $(-1,0)$.

Let $\tilde t$ be such that $(u,v)(t)$ is close to $(-1,0)$ for $t\geq \tilde t$. Then, $\frac {v(t)}{1+u(t)}\geq\min\left( \frac {v(\tilde t)}{1+u(\tilde t)},\frac {-\sqrt 2}{3B}\right)$ for all $t\geq \tilde t$, and thus, $u'(t)=v(t)\geq -C(u(t)+1)$, which implies the estimate
\begin{equation*}
u(t)\geq -1+(u(\tilde t)+1)e^{-C(t-\tilde t)}.
\end{equation*}
Since $(u,v)(t)\to (-1,0)$ as $t\to \bar t$, it follows that $\bar t=+\infty$, and $(u,v)$ is indeed globally defined.


\bigskip

\emph{Step 2:} We show next that there exists a constant $\nu_B$ such that, there exists a propagative front if and only if $\nu>\nu_B$.

If $\nu>\frac {\sqrt 2} {B}$, we have $F_V(U,-B/\sqrt 2)=\nu B/\sqrt 2-U> 0$ for $U\in (-1,1)$. Therefore, the solution given by Lem. \ref{l1} cannot cross the line $V=-B/\sqrt 2$, and thus it defines a propagative front thanks to Step 1. We deduce that, there exists a propagative front if $\nu$ is large enough.

\medskip

On $[-1/\sqrt{2},0]\times [-B/\sqrt 2,0]$, we have, for $\nu\leq -(1+2\sqrt 2)B$,
\begin{equation*}
F_V(U,V)\leq -\nu V-4V(V+B/ \sqrt 2)\leq BV\leq B\,F_U(U,V).
\end{equation*}
It follows that the solution given by Lem. \ref{l1} necessarily crosses the line $V=-B$, and thus it does not define a propagative front thanks to Step 1. We deduce that the model does not admit a propagative front if $-\nu$ is large enough.

\medskip

Consider a solution $(u_{\nu_1},v_{\nu_1})(t)$ given by Lem. \ref{l1}  for some $\nu_1$, that converges to $(-1,0)$ as $t\to\infty$, and $\nu_2> \nu_1$. Then, since $\frac {\sqrt{\nu^2+4}-\nu}2$ is a decreasing function of $\nu$ and 
\begin{equation*}
(u_{\nu_i},v_{\nu_i})(t)\sim-C_-e^{\frac {\sqrt{\nu_{\nu_i}^2+4}-\nu_{\nu_i}}2}\left(1,\frac {\sqrt{\nu_{\nu_i}^2+4}-\nu_{\nu_i}}2\right),
\end{equation*}
the graph of $(u_{\nu_1},v_{\nu_1})(t)$ is below the graph of $(u_{\nu_2},v_{\nu_2})(t)$ for $t<<0$. Thanks to Lem. 2, this implies that the whole graph of $(u_{\nu_1},v_{\nu_1})$ is below the graph of $(u_{\nu_2},v_{\nu_2})$. Using the latter and Step 1 we obtain that $(u_{\nu_2},v_{\nu_2})$ defines a propagative front.
%

\medskip

Finally, we show that $\nu_B$ is a decreasing function of $B$. Firstly we notice that, for $B_1\leq B_2$, we have $F_U^{B_1}=F_U^{B_2}$, and $F_U^{B_1}\leq F_U^{B_2}$ on $(-1,0]\times \mathbb R_-$. It follows that, thanks to Step 1, if the solution given by Lem. \ref{l1} for $B_1$ converges to $(-1,0)$ as $t\to\infty$, so does the one associated to $B_2$. This shows that $\nu_B$ is a decreasing function of $B$.

\begin{flushright}
$\square$
\end{flushright}

\bigskip

\noindent\textbf{Proof of Prop. \ref{EDOstat}:}

Assume that $\nu_B<0$. Then, for $\nu=\frac {\nu_B}2$, the solution $(u_{\nu_B/2},v_{\nu_B/2})$  given by Lem. \ref{l1}, is globally defined, and it satisfies $(u_{\nu_B/2},v_{\nu_B/2})(t)\to (-1,0)$ as $t\to +\infty$. Moreover we have
\begin{equation*}
(u_{\nu_B/2},v_{\nu_B/2})(t)\sim_{t\to -\infty} - C_- e^{\frac {\sqrt{(\nu_B/2)^2+4}-{\nu_B/2}}2t}\left(1,\frac {\sqrt{(\nu_B/2)^2+4}-{\nu_B/2}}2\right),
\end{equation*}
and $\frac {\sqrt{(\nu_B/2)^2+4}-{\nu_B/2}}2>1$. 

\medskip

Consider now the vector field \eqref{FUV} for $\nu=0$. Since $(u_0,v_0)(t)\sim_{t\to -\infty} - C_-' e^{t}(1,1)$, for $\bar t$ small enough, $(u_0,v_0)(\bar t)$ is strictly above the graph of $(u_{\nu_B/2},v_{\nu_B/2})$. Let $\tilde v$ be such that $(u_0(\bar t),\tilde v)$ is strictly between $(u_0,v_0)(\bar t)$ and the graph of $(u_{\nu_B/2},v_{\nu_B/2})$. We define $(\bar u,\bar v)$ to be the solution of the ODE given by the vector fields \eqref{FUV} such that $(\bar u,\bar v)(0)=(u_0(\bar t),\tilde v)$ and $\nu=0$. Then, thanks to Lem. \ref{l2}, $(\bar u,\bar v)$ is defined on $\mathbb R_+$. Moreover, since  $(u_0,v_0)(t)\to (0,0)$ as $t\to -\infty$ and $(0,0)$ is a hyperbolic point, there exists $\tilde t>0$ such that $\bar u(\tilde t)=0$. 

By symmetry, $(\bar u,\bar v)(\tilde t+t)=(-\bar u,\bar v)(\tilde t-t)$, and thus, $(\bar u,\bar v)$ is globally defined and satisfies $(\bar u,\bar v)(t)\to (\pm 1,0)$ as $t\to \pm \infty$. This completes the proof of proposition \ref{EDOstat}.

\begin{flushright}
$\square$
\end{flushright}

\section{The case of asexual populations}\label{section:as}

\subsection{The model}

We consider here the same fitness \eqref{fitness} as we considered for sexual populations. If we additionally model mutations through a diffusion of rate $\sigma_v^2$ (for more on the different ways to model mutations, see \cite{CFM}), then the evolution of a population structured by a phenotypic trait $v$ only can be modeled by the classical model (see e.g. \cite{Burger, DBLD}):
\begin{equation*}
\partial_t n(t,v)-\sigma_v^2\Delta_vn(t,v)=\left[r_{max}-\frac 1{2V_s}(v-\theta)^2-\frac 1K\int n(t,w)\,dw\right]n(t,v).
\end{equation*}
Just as in the sexual case, we add a spatial structure to this model, the population is then structured by both a phenotypic trait $v$ as above, but also by a space variable $x\in \mathbb R$: $n(t,x,v)$. We assume that the selection-mutation process described above occurs locally in space, but that individuals move randomly in space (which we model by a diffusion of rate $\sigma_x$), and that the trait of optimal adaptation changes linearly in space (see \eqref{b}). We then get the following model for asexual populations ($\frac {I_s}2$ represents lethal mutations):
\begin{align*}
&\partial_t n(t,x,v)-\sigma_x^2\Delta_xn(t,x,v)-\sigma_v^2\Delta_vn(t,x,v)\\
&\quad=\left[r_{max}-\frac {I_s}2-\frac 1{2V_s}(v-bx)^2-\frac 1K\int n(t,x,w)\,dw\right]n(t,x,v).
\end{align*}
We then rescale the problem as follows:
\begin{equation*}
\tilde t=\frac 1 {r_{max}-\frac {I_s}2}t, \quad \tilde n=K(r_{max}-\frac {I_s}2)n,
\end{equation*}
\begin{equation*}
\tilde x=\frac{\sqrt{r_{max}-\frac {I_s}2}}{\sigma_x}x,\quad \tilde v=\frac{\sqrt{r_{max}-\frac {I_s}2}}{\sigma_v}v,
\end{equation*}
and define the two parameters
\begin{equation*}
A:=\frac 1{2V_s\sigma_v\sqrt{r_{max}-\frac {I_s}2}},\quad B:=\frac {\sigma_v}{\sigma_x}b.
\end{equation*}
Then, we obtain the following rescaled model:
\begin{align}
&\partial_tn(t,x,v)-\Delta_xn(t,x,v)-\Delta_vn(t,x,v)\nonumber\\
&\quad =\left[1-A(v-Bx)^2-\int n(t,x,v')\,dv'\right]n(t,x,v)\label{eqa}
\end{align}

For an existence theory for this equation, we refer to \cite{ADP,Prevost}. Notice that an integration of \eqref{eqa} over the $v$ variable provides the uniform in time estimate:
\begin{equation}\label{Na}
\|n\|_{L^\infty(t,x,L^1(v))}\leq \max\left(\|n^0\|_{L^\infty(x,L^1(v))},\,1\right).
\end{equation}

\subsection{Qualitative properties of asexual population}


We show that if $A(1+B^2)>1$, then, for any initial population, the population goes extinct when $t\to \infty$:
\begin{prop}\label{exta}
Assume $A(1+B^2)>1$. For any initial population $n^0\in L^\infty$, the population will go extinct exponentially fast when $t\to\infty$:
\begin{equation*}
\left\|e^{\frac{1}2\sqrt{\frac A {1+B^2}}(v-Bx)^2}n(t,\cdot,\cdot)\right\|_{L^\infty(\mathbb R^2)}=O(e^{-ct}),
\end{equation*}
where $c=\sqrt{A(1+B^2)}-1$.
\end{prop}

If $A(1+B^2)<1$, the population survives, and does not remain confined in a given part of the space set:
\begin{prop}\label{extb}
Assume $A(1+B^2)<1$. There exists $X,\,\kappa>0$ depending only on $A$ and $B$ such that for any initial population $n^0\neq 0$, any $T_1>0,\,x_0\in\mathbb R$, there exists $T_2\geq T_1$ such that  
\begin{equation*}
\|n(T_2,\cdot,\cdot)\|_{L^1([x_0-X,x_0+X]\times \mathbb R)}>\kappa.
\end{equation*}
\end{prop}

\begin{rem}
We believe that if $A(1+B^2)<1$, then the population indeed invades the whole space in the stronger sense that there exists $X,\,\kappa>0$ depending only on $A$ and $B$ such that for any initial population $f^0\neq 0$, and any $x_0\in\mathbb R$, there exists $T\geq 0$ such that for all $t\geq T$,
\begin{equation*}
\|n(t,\cdot,\cdot)\|_{L^1([x_0-X,x_0+X]\times \mathbb R)}>\kappa.
\end{equation*}
We were unfortunatly unable to show this stronger result.

On the contrary, it is not clear that the result of Prop. \ref{extb} would be true locally, that is if $X$ could be choosen arbitrarly small. It is indeed known that reaction-diffusion equations with integral terms may lead to complicated invasion fronts , see \cite{BNPR}.

\end{rem}

\begin{rem}\label{eq2as}
Prop. \ref{exta} and Prop. \ref{extb} show that the dynamics of asexual populations is different from the one of sexual populations: For sexual populations, numerical simulations from \cite{KB} show that for some parameters, a population can survive, but remain confined in a restricted area (this observation being comforted by Prop. \ref{EDOstat}, although the model \eqref{eq:modelZ} is too simple to model extinction cases). This shows a property that is well known experimentally: asexual populations have a wider geographic distributions than asexual populations (see \cite{PYW,BESP}).

Surprisingly, Prop. \ref{propNZ} also holds for asexual populations (see also Rem. \ref{geneflow}). The reason why the model \eqref{eq2} cannot be used to model asexual population is that the closure condition \eqref{closure} is based on the Assumption 1, that is on the sexual reproduction kernel, which cannot be satisfied in the case of asexual populations.
\end{rem}

\bigskip

\noindent\textbf{Proof of Prop. \ref{exta}:}

\medskip

We consider the following function, for some $\lambda>0$:
\begin{equation*}
\phi(t,x,v):=\lambda e^{-(\sqrt{A(1+B^2)}-1)t}e^{-\frac 12\sqrt{\frac A{1+B^2}}(v-Bx)^2}
\end{equation*}
that is a strong solution of:
\begin{equation}
\partial_t\phi(t,x,v)-\Delta_x\phi(t,x,v)-\Delta_v\phi(t,x,v) =\left[1-A(v-Bx)^2\right]\phi(t,x,v).\label{maj1}
\end{equation}
Notice that since $-n(t,x,v)\,\int n(t,x,v')\,dv'\leq 0$, $n$ is a subsolution of \eqref{maj1}. Moreover if we choose $\lambda:=\left\|e^{\frac{1}2\sqrt{\frac A {1+B^2}}(v-Bx)^2}n^0\right\|_\infty$, then $n^0\leq \phi(0,\cdot,\cdot)$. Thanks to the comparison principle for the equation \eqref{maj1}, we obtain
\begin{equation*}
n(t,x,v)\leq \phi(t,x,v),\quad \forall t\geq 0,\,x\in\mathbb R,\,v\in\mathbb R.
\end{equation*}

Since $A(1+B^2)>1$, $\phi$ vanishes as $t$ tends to infinity:
\begin{eqnarray*}
\|e^{\frac{1}2\sqrt{\frac A {1+B^2}}(v-Bx)^2}n(t,\cdot,\cdot)\|_{L^\infty(x,L^1(v))}&\leq&\|e^{\frac{1}2\sqrt{\frac A {1+B^2}}(v-Bx)^2}\phi(t,\cdot,\cdot)\|_{L^\infty(x,L^1(v))}\\
&\leq& C\,e^{-(\sqrt{A(1+B^2)}-1)t}\to 0,\textrm{ as }t\to\infty.
\end{eqnarray*}
This concludes the proof of Prop. \ref{exta}.
\begin{flushright}
$\square$
\end{flushright}

\bigskip

\noindent\textbf{Proof of Prop. \ref{extb}:}

\medskip

\emph{Step 1:} We show that $\|n(\bar t+1,\cdot,\cdot)\|_{L^\infty([x_0-X/2,x_0+X/2],L^1(v))}$ can be controlled by $\|n(\bar t,\cdot,\cdot)\|_{L^1([x_0-X,x_0+X]\times \mathbb R}$ for some $X>0$.

To show this, we notice that $\varphi(t,x,v):=\frac{e^t}{4\pi t}e^{\frac{-(x^2+v^2)}{4t}}$ is the fundamental solution of:
\begin{equation}
\partial_t\varphi(t,x,v)-\Delta_x\varphi(t,x,v)-\Delta_v\varphi(t,x,v) =\varphi(t,x,v).\label{maj2},
\end{equation}
and that $n$ is a subsolution of the same equation. Then, the comparison principle for \eqref{maj2} shows that $$n(\bar t+1, \cdot,\cdot)\leq(\varphi(1)\ast_{x,v}n(\bar t))(\cdot,\cdot),$$
since $n(\bar t,\cdot,\cdot)\leq (\varphi(\bar t)\ast_{x,v}n^0)(\cdot,\cdot)$. In particular, we have
\begin{align*}
&\|n(\bar t+1,\cdot,\cdot)\|_{L^\infty([x_0-X/2,x_0+X/2],L^1(v))}\\
&\quad\leq\|\varphi(1)\ast_{x,v}n(\bar t\,)\|_{L^\infty([x_0-X/2,x_0+X/2],L^1(v))}\\
&\quad\leq\|\varphi(1)\ast_{x,v}(n(\bar t\,)\mathds{1}_{[x_0-X,x_0+X]\times \mathbb R })\|_{L^\infty([x_0-X/2,x_0+X/2],L^1(v))}\\
&\qquad+\|\varphi(1)\ast_{x,v}(n(\bar t\,)\mathds{1}_{[x_0-X,x_0+X]^c\times \mathbb R })\|_{L^\infty([x_0-X/2,x_0+X/2],L^1(v))}\\
&\quad\leq \|n(\bar t\,)\|_{L^1([x_0-X,x_0+X]\times \mathbb R)}\,\|\varphi(1,\cdot,\cdot)\|_{L^\infty(x,L^1(v))}\\
&\qquad+\|n(\bar t\,)\|_{L^\infty(x,L^1(v))}\int_{|x-x_0|\geq X/2}\int\varphi(1,x,v)\,dv\,dx\\
&\quad\leq\frac e{\sqrt{4\pi}} \|n(\bar t\,)\|_{L^1([x_0-X,x_0+X]\times \mathbb R)}+\|n(\bar t\,)\|_{L^\infty(x,L^1(v))}\frac{2e}{\sqrt{4\pi}}\int_{X/2}^\infty e^{\frac{-x^2}4}\,dx\\
&\quad\leq C \|n(\bar t\,)\|_{L^1([x_0-X,x_0+X]\times \mathbb R)}+o_{X\to\infty}(1)\|n(\bar t\,)\|_{L^\infty(x,L^1(v))}.
\end{align*}

\bigskip

\emph{Step 2:}  We show that if $t\geq 0$ is large enough, then $n$ can be minored by a Gaussian function.

We define:
\begin{equation*}
\psi(t,x,v):=\frac{e^{-ct}}{4\pi t}e^{-\frac{x^2+v^2}{4t}}e^{-\alpha\frac{(v-Bx)^2}2}.
\end{equation*}

Then,
\begin{align*}
&(\partial_t-\Delta_x-\Delta_v)\psi(t,x,v)\\
&\quad=\left[\left(\alpha(1+B^2)-c\right)-\alpha^2(1+B^2)(v-Bx)^2-\alpha\frac {(v-Bx)^2}{2t}\right]\psi(t,x,v)\\
&\quad\leq\left[\left(\alpha(1+B^2)-c\right)-\alpha^2(1+B^2)(v-Bx)^2\right]\psi(t,x,v).
\end{align*}

Let $x_0\in\mathbb R$, and $R>0$ large enough for $\int n^0_R(x,v)\,dx\,dv>0$ to hold, where
\begin{equation*}
n^0_R=n^0|_{\{(x,v);|x-x_0|^2+|v|^2\leq R^2\}}.
\end{equation*}
We define $\tilde n(t,x,v):=\left(\psi(t,\cdot,\cdot)\ast_{x,v}n^0_R\right)(x,v)$. Then,
\begin{align*}
&(\partial_t-\Delta_x-\Delta_v)\tilde n(t,x,v)\\
&\quad=((\partial_t\psi(t,\cdot,\cdot)-\Delta_x\psi(t,\cdot,\cdot)-\Delta_v\psi(t,\cdot,\cdot))\ast_{x,v}n^0_R)(t,x,v)\\
&\quad\leq\left(\left[\left(\left(\alpha(1+B^2)-c\right)-\alpha^2(1+B^2)(v-Bx)^2\right)\psi(t,\cdot,\cdot)\right]\ast_{x,v}n^0_R\right)(t,x,v)\\
&\quad\leq \left(\alpha(1+B^2)-c\right)\tilde n(t,x,v)-\alpha^2(1+B^2)\left(((v-Bx)^2\psi(t,\cdot,\cdot))\ast_{x,v}n^0_R\right)(t,x,v).
\end{align*}
we can estimate the last term using the fact that $\textrm{supp }n^0_R$ is bounded, and a Young inequality:
\begin{align*}
&-\left(((v-Bx)^2\psi(t,\cdot,\cdot))\ast_{x,v}n^0_R\right)(t,x,v)\\
&\quad =-\int \int\left((v-Bx)-(v'-Bx')\right)^2\psi(t,x-x',v-v')n^0_R(x',v')\,dv'dx'\\
&\quad =-\int \int\left((v-Bx)^2+(v'-Bx')^2-2(v-Bx)(v'-Bx')\right)\\
&\phantom{\quad =-\int \int dsf}\psi(t,x-x',v-v')n^0_R(x',v')\,dv'dx'\\
&\quad \leq \left[-(v-Bx)^2+|v-Bx|\right]\tilde n(t,x,v)\\
&\quad \leq \left[-(1-\delta)(v-Bx)^2+ \frac 1\delta\right]\tilde n(t,x,v).
\end{align*}
Finally, we get that:
\begin{align*}
&(\partial_t-\Delta_x-\Delta_v)\tilde n(t,x,v)\\
&\leq \left[\alpha(1+B^2)(1+\frac \alpha\delta)-c-(1-\delta)\alpha^2(1+B^2)(v-Bx)^2\right]\tilde n(t,x,v),
\end{align*}
We now choose $\alpha$ and $c$ as follows:
\begin{equation*}
\alpha:=\sqrt{\frac {A}{(1+B^2)(1-\delta)}}
\end{equation*}
\begin{eqnarray*}
c&=&\alpha(1+B^2)(1+\frac \alpha\delta)+ \max\left(\|n^0\|_{L^\infty(x,L^1(v))},\,1\right)\\
&=&\sqrt{\frac{A(1+B^2)}{1-\delta}}(1+B^2)+\frac {A}{\delta(1-\delta)}+ \max\left(\|n^0\|_{L^\infty(x,L^1(v))},\,1\right),
\end{eqnarray*}
then $\tilde n$ is a subsolution of:

\begin{align}
&\partial_t\psi(t,x,v)-\Delta_x\psi(t,x,v)-\Delta_v\psi(t,x,v)\nonumber \\
&\quad =\left[0-A(v-Bx)^2- \max\left(\|n^0\|_{L^\infty(x,L^1(v))},\,1\right)\right]\psi(t,x,v).\label{min1}
\end{align}
Since $\|n\|_{L^\infty(t,x,L^1(v))}\leq \max\left(\|n^0\|_{L^\infty(x,L^1(v))},\,1\right)$, $n$ is a supersolution of \eqref{min1}, and since $n^0\geq  f^0_R$, the comparison principle for \eqref{min1} shows that $n\geq\tilde n$, for all $t\geq 0$, and in particular,
\begin{eqnarray*}
n(t,x,v)&\geq&\min_{\|(y,w)\|\leq R}\psi(t,(x-x_0)+y,v+w)\,\int_{\|(y,w)\|\leq R}n^0_R(x_0+y,w)\,dy\,dw\\
&\geq&C\frac{e^{-ct}}{4\pi t}\min_{\|(y,w)\|\leq R}\left(e^{-\frac{(x-x_0+y)^2+(v+w)^2}{4t}}e^{-\alpha\frac{(v+w-B(x-x_0+y))^2}2}\right)\\
&\geq&C\frac{e^{-ct}}{4\pi t} \exp\Big\{-\frac 1{4t}\big((x-x_0)^2+R^2+2R|x-x_0|+(v-Bx)^2\\
&&\qquad+2(B|x-x_0|+B|x_0|+R)|v-Bx|+B^2(x-x_0)^2\\
&&\qquad+2B(B|x_0|+R)|x-x_0|+(B|x_0|+R)^2\big)\Big\}\\
&&\exp\left\{-\alpha\left((v-Bx)^2+2(R+BR+B|x_0|)\,|v-Bx|+(R+BR+B|x_0|)^2\right)\right\},
\end{eqnarray*}
and then, thanks to Young inequalities, we get, for any $\delta'>0$:
\begin{eqnarray*}
n(t,x,v)&\geq& C_{\delta'}\frac{e^{-ct}}{4\pi t} \exp\left\{-\frac {C}{4t}(x-x_0)^2\right\}\\
&&\exp\left\{-\left(\alpha+\delta'+\frac C{4t}\right)(v-Bx)^2\right\},
\end{eqnarray*}
where $C_{\delta'}$ depends on $\delta',\,x_0,\,B,\,\int_{\|(y,w)\|\leq R}n^0(x_0+y,w)>0$. Finally, for any $x_0\in\mathbb R$, $\mu_1,\mu_2>0$ and $T_1\geq 0$, there exists $\tau>T_1,\,\lambda=\lambda(\tau)>0$ such that
\begin{equation*}
n(\tau,x,v)\geq \lambda e^{-\frac 12\left(\sqrt{\frac A{1+B^2}}+\mu_1\right)(v-Bx)^2}e^{-\mu_2 \frac{(x-x_0)^2}2},\quad \forall x,\,v\in \mathbb R\times \mathbb R.
\end{equation*}

\bigskip

\emph{Step 3:} We show the result.

Since $A(1+B^2)<1$, there exists $\mu>0$ such that (this defines the $\mu$ and then the $\kappa$):
\begin{equation}\label{signec}
1-\sqrt{A(1+B^2)}-\mu(3+B^2)>0.
\end{equation}
We define:
\begin{equation*}
\chi(t,x,v):=e^{\left(1-\sqrt{A(1+B^2)}-\mu(3+B^2)\right)t}e^{-\frac 12\left(\sqrt{\frac A{1+B^2}}+\mu\right)(v-Bx)^2}e^{-\frac\mu 2x^2},
\end{equation*}
which satisfies:
\begin{align*}
&\partial_t\chi(t,x,v)-\Delta_x\chi(t,x,v)-\Delta_v\chi(t,x,v)\\
&\quad=\left[1-\mu-(\sqrt A+\mu\sqrt{1+B^2})^2(v-Bx)^2-\mu^2 x^2\right]\chi(t,x,v)\\
&\quad\leq\left[1-\mu-A(v-Bx)^2-\max\left(\|n^0\|_{L^\infty(x,L^1(v))},\,1\right)1_{|x-x_0|>X}\right]\chi(t,x,v),
\end{align*}
if $X$ is large enough. Then, $\chi$ is a subsolution of \eqref{min2}:
\begin{align}
&\partial_t\chi(t,x,v)-\Delta_x\chi(t,x,v)-\Delta_v\chi(t,x,v) \nonumber\\
&\quad =\left[1-\mu-A(x-Bv)^2-\max\left(\|n^0\|_{L^\infty(x,L^1(v))},\,1\right)1_{|x-x_0|>X}\right]\chi(t,x,v).\label{min2}
\end{align}
 
\medskip

Thanks to Step 1, provided that $X$ is large, there exists $\kappa>0$ such that if $\|n(t,\cdot,\cdot)\|_{L^1([x_0-X,x_0+X]\times \mathbb R)}\leq \kappa$ for all times $t\geq T_1$, then, for any $t\geq T_1+1$, $n$ is bounded as follows:
\begin{equation*}
\left\|\int n(t,\cdot,v)\,dv\right\|_{L^\infty([x_0-X/2,x_0+X/2])}\leq \mu.
\end{equation*}
Then, $t\mapsto n(T_1+1+t,\cdot, \cdot)$ is a super-solution of \eqref{min2} for $t\geq 0$.



Thanks to step 2, there exists $\tau>T_1$ and $\lambda>0$ such that 
\begin{equation*}
n(\tau,x,v)\geq \lambda\, \chi(0,x,v),\quad \forall x,\,v\in \mathbb R\times \mathbb R.
\end{equation*}
Then, thanks to the comparison principle, $n(\tau + t,\cdot,\cdot)\geq \lambda \chi(t,\cdot,\cdot)$ for $t\geq 0$, and in particular,
\begin{eqnarray*}
\|n(\tau+t,\cdot,\cdot)\|_{L^1}&\geq& \lambda \|\chi(t,\cdot,\cdot)\|_{L^1}\\
&\geq& C\, e^{\left(1-\sqrt{A(1+B^2)}-\mu(3+B^2)\right)t}\to\infty,
\end{eqnarray*}
thanks to \eqref{signec}. This shows that  the assertion that $\|n(t,\cdot,\cdot)\|_{L^1([x_0-X,x_0+X]\times \mathbb R)}\leq \kappa$ for all times $t\geq T_1$ cannot be valid, and shows Prop. \ref{extb}.
\begin{flushright}
$\square$
\end{flushright}

\bigskip

\textbf{Acknowledgments:} GR has been supported by Award No. KUK-I1-007-43 of Peter A. Markowich, made by King Abdullah University of Science and Technology (KAUST).

\end{document}